\documentclass[11pt,oneside]{amsart}
\usepackage{amssymb,amsmath,latexsym,amscd,amsfonts}
\usepackage{graphics}
\usepackage{epsfig}
\usepackage{psfrag}
\def\ignore #1 {}

\newtheorem*{ack}{Acknowledgments}
\def\hpic #1 #2 {\mbox{$\begin{array}[c]{l} \epsfig{file=#1,height=#2} \end{arr\
ay}$}}
\def\vpic #1 #2 {\mbox{$\begin{array}[c]{l} \epsfig{file=#1,width=#2} \end{arra\
y}$}}

\newcommand{\yco}[2]{Y^{[#1]}_{(#2)}}

\def\M{\mathbf M}
\def\Z{\mathbb Z}
\def\C{\mathbb C}
\def\R{\mathbb R}
\def\P{\mathbb P}

\def\RR{\mathcal R}
\def\ZZ{\mathbf Z}

\begin{document}


\title{Twisted homogeneous coordinate rings of abelian surfaces via mirror symmetry}

\author{Marco Aldi}
\begin{abstract} 
In this paper we study Seidel's mirror map for abelian and Kummer surfaces. We find that mirror symmetry leads in a very natural way to the classical parametrization of Kummer surfaces in $\P^3$. Moreover, we describe a family of embeddings of a given abelian surface into noncommutative projective spaces.
\end{abstract}
\maketitle

\section{Introduction}
\label{intro}

Let $Y$ be a symplectic Calabi-Yau manifold with symplectic form $\omega$, $L\subset Y$ a Lagrangian submanifold and and $\rho$ a symplectomorphism. Using these data, Seidel proposed a method to explicitly construct a complex Calabi-Yau manifold $X$ mirror to $Y$, provided that $X$ is projective. In fact, assuming the homological mirror symmetry conjecture \cite{K}, we could write the homogeneous coordinate ring of $X$ in terms of the derived Fukaya category of $Y$ (see \cite{Z} for more details) as 
$$
\RR(Y) := \bigoplus_{k\ge 0} {\rm Hom}_{\rm DFuk(Y)} (L,\rho^k L) \,.
$$
Even if homological mirror symmetry is not proven for $Y$ we can {\it define} $X$ as the projective spectrum of $\RR(Y)$. 

Seidel's construction has been investigated for elliptic curves \cite{Z}, abelian varieties of higher dimension \cite{AZ} and toric Fano varieties \cite{A}. 

In the first part of this paper we study Seidel's mirror map for a symplectic Kummer surface $K$, i.e. a quotient of symplectic four-torus with respect to the involution reversing the sign on each coordinate. Using the methods of \cite{AZ}, we are able to show that $\RR(K)$ is the homogeneous coordinate ring of a singular complex surface of Kummer type embedded in $\P^3$. The mirror correspondence is made very explicit as we are able to express the coefficients of the quartic polynomial defining the complex surface in terms of theta functions (on abelian surfaces) depending on the original symplectic form. Incidentally, we find a natural symplectic interpretation for the classical (see \cite{H}, \cite{GD}) parametrization of the space of Kummer surfaces in $\P^3$ in terms of $(16,6)$-configurations in $\P^3$.

In \cite{AZ} it was shown that Seidel's mirror map can be used to reconstruct twisted homogeneous coordinate rings or, equivalently, embeddings of the complex mirror into an ambient noncommutative projective variety \cite{SVdB}. For example, if $Y=E$ is a symplectic two-torus and $\rho$ is affine linear in the universal cover, then $\RR(E)$ turns out to be isomorphic to a quotient of the 3-dimensional Sklyanin algebra \cite{ATV}, viewed as the twisted homogeneous coordinate ring of noncommutative $\P^2$. 

In the second part of this paper we consider the case where $Y=T$ is a symplectic four-torus and $\rho$ is affine linear in the universal cover. The corresponding twisted homogeneous coordinate ring is a quotient of a noncommutative deformation of the (commutative) homogenous polynomial ring in 9 variables. This family of deformations is parametrized by a choice of symplectic form on the four-torus $T$ and a point on $T$.  We interpret these noncommutative rings as twisted homogeneous coordinate rings of a family of noncommutative deformations of $\P^8$ in which the mirror of $T$ is embedded. 

The role of elliptic curves in the original definition of the 4-dimensional Sklyanin algebra \cite{Sk} was due to Sklyanin's interest in the classification of quantum integrable systems constructed using Baxter's solution to the quantum Yang-Baxter equation. Baxter's solution is parametrized by elliptic functions, a fact geometrically explained by Cherednik \cite{C}. Moreover, in \cite{C} more general solutions to the quantum Yang-Baxter equation were found in terms of theta functions on higher dimensional abelian varieties. We consider it natural to look at the family of noncommutative deformations of $\P^8$ found in this paper as an analogue of Sklyanin algebras based on Cherednik's $R$-matrices. It would be of interest to study homological properties of these twisted homogenous coordinate rings in the spirit of \cite{ATV}, \cite{S}, \cite{SVdB}.

\begin{ack}
This work was initiated while preparing \cite{AZ}. I would like to thank Prof. Eric Zaslow for many helpful discussions and constant encouragement. This work was partially supported by NSF grant DMS-0072508.    
\end{ack}

\section{Seidel's mirror map for Kummer surfaces in $\P^3$}\label{Kummer}

\subsection{Kummer surfaces in $\P^3$} \label{kummer}
Here we recall some classical facts (\cite{H}, see \cite{GD} for a modern and rigorous treatment).
By definition, a (singular) Kummer surface is the quotient of an abelian surface with respect to the standard involution reversing the sign of each of the coordinates. Kummer surfaces embedded in $\P^3$ bijectively correspond  to singular quartic hypersurfaces whose only singularities are sixteen double points. It is convenient to parametrize Kummer surfaces in $\P^3$ by points $(g:h:j:k)\in \P^3$ such that
$$
\displaylines{
gh \neq \pm jk, \quad gj \neq \pm hk , \quad gk \neq \pm hj\, , \cr
g^2 + h^2 \neq j^2 + k^2, \quad g^2 + j^2 \neq h^2 + k^2, \quad g^2 + k^2 \neq h^2 + j^2\, , \cr
g^2 + h^2 + j^2 + k^2  \neq  0\, . 
}
$$
In fact, it is possible to choose coordinates such that every Kummer surface in $\P^3$ has the equation
$$
\sum_{i=0}^3 X_i^4 + A(X_0^2 X_1^2 + X_2^2 X_3^2) + B(X_0^2 X_2^2 + X_1^2 X_3^2) + C (X_0^2 X_3^2 + X_1^2 X_2^2) + D X_0 X_1 X_2 X_3 = 0
$$
where 
$$
\displaylines{
A = \frac{j^4 + k^4 - h^4 - g^4}{h^2 g^2 - j^2 k^2}, \quad B = \frac{k^4 + h^4 - j^4 - g^4}{j^2 g^2 - h^2 k^2}, \quad C = \frac{h^4 + j^4 - k^4 - g^4}{k^2 g^2 - j^2 h^2}\, ,\cr
D = 2\frac{hjkg(g^2 + h^2 - j^2 - k^2)(g^2+j^2-k^2-h^2)(g^2-h^2-j^2+k^2)(h^2+j^2+k^2+g^2)}{(h^2g^2-j^2k^2)(j^2g^2-h^2k^2)(k^2g^2-j^2h^2)}\, . 
}
$$

\subsection{Seidel's Mirror Map}\label{Seidel}

Let $T:=\R^4/\Z^4$ be a real four-torus endowed with complex symplectic form 
\begin{equation}\label{omega}
\omega=\tau_1 dx_1 \wedge dy_1 + \tau_2 dx_2 \wedge dy_2 + \tau_3 (dx_1\wedge dy_2 + dx_2 \wedge dy_1)\, ,
\end{equation}
and consider the standard involution $\iota(x_1, y_1, x_2, y_2) = - (x_1, y_1, x_2, y_2)$. We call the quotient symplectic manifold $K:= T/\iota$ a symplectic Kummer surface.  We compute the homogenous coordinate ring $\RR(K)$
of the mirror by choosing the symplectomorphism to be
$$
\rho(x_1,y_1,x_2,y_2) = (x_1,x_2, y_1 + 2 x_1, y_2 + 2 x_2)
$$
and the Lagrangian $L\subset K$ that can be lifted to 
$$
\widetilde L:=\{y_1 = 0 , y_2 = 0 \}
$$
as base of the brane fibration in the universal cover $\R^4$. The points
$$
\yco{l}{a,b} : = (a/2l,0,b/2l,0)\in L\cap \rho^l L
$$
form a linear basis for $\RR(T)$, where $l\in \Z_{\ge 0}$ and $a,b\in \Z/2l\Z$. The basic product formula for $\RR(T)$ (see \cite{AZ}) is
\begin{equation}\label{prod}
\yco{l}{a,b}\yco{l}{c,d}:=\sum_{i, j \in \Z/2\Z \times \Z/2\Z}  \theta(c - a + 2 i l, d - b + 2 j l, l) \yco{2 l}{a + c + 2 l i, b + d + 2 l j}
\end{equation}
where
$$
\theta(a,b,l):= \sum_{m,n\in \Z} e^{4 l \pi i \tau_1 (\frac{a}{4l} + m)^2} e^{4 l \pi i \tau_2 (\frac{b}{4 l} + n)^2} e^{8 l\pi i \tau_3(\frac{a}{4l} + m) (\frac{b}{4l} + n)}\, . 
$$
Since $K$ is an orbifold, we need to explain how to account for the group action in the definition of $\RR(K)$. Let $\phi:\RR(T)\to\RR(K)$ be the map graded linear spaces that identifies generators in the same orbit under the action of $\iota$. It is easy to see that in degree $l$ there are $2(l^2+1)$ distinct images $\phi(\yco{l}{a,b})$. We use $\phi$ to induce a product on $\RR(K)$ according to the formula
$$
\phi(Y_1)\phi(Y_2)=\phi[(\sum_{Y'_1\in \phi^{-1}Y_1} Y'_1)(\sum_{Y'_2\in \phi^{-1} Y_2} Y'_2)]
$$ 
where $Y_1$ and $Y_2$ are generators of $\RR(K)$. In particular, this means that constant holomorphic disks located at the fixed points of $\iota$ contribute by $1/2$ to the generating function for the structure constants of products in $\RR(K)$. For the rest of the section we work exclusively in $\RR(K)$ and so we simplify the notation by omitting $\phi$.

We denote by
$$
X_0:= \yco{1}{0,0}, \quad X_1:= \yco{1}{1,0}, \quad X_2:= \yco{1}{0,1} , \quad X_3:=\yco{1}{1,1}
$$
the generators of $\RR(K)$ in degree one. In degree $l$ there are ${3+l \choose l}$ homogeneous monomials in four variables, and since ${3+l \choose l} - 2(l^2+1)={l-1\choose l-4}$, we expect only one nontrivial relation in degree $l=4$.

Consider the following five quartic polynomials \footnote{This choice can be justified {\it a priori}. In fact by looking at which generators of $\RR(K)$ contribute to a given quartic monomial, on can see that any other relation involving other quartic monomials would imply a (non existing) nontrivial relation in degree two.} 
$$
\displaylines{
W_0 := X_0^4 + X_1^4 + X_2^4 + X_3^4, \quad W_4:=X_0 X_1 X_2 X_3\, , \cr
W_1 := X_0^2 X_1^2 + X_2^2 X_3^2, \quad W_2 := X_0^2 X_2^2 + X_1^2 X_3^2, \quad W_3 := X_0^2 X_3^2 + X_1^2 X_2^2 
}
$$
and the four linear combinations
$$
\displaylines{
Z_0 := \yco{4}{0,0} + \yco{4}{4,0} + \yco{4}{0,4} + \yco{4}{4,4}\, , \cr
Z_1 :=  \yco{4}{2,0} + \yco{4}{2,4}, \quad Z_2 :=  \yco{4}{0,2} + \yco{4}{4,2} , \quad Z_3 :=  \yco{4}{2,2} + \yco{4}{2,6}  
}
$$
of elements of $\RR(K)$ in degree four. As a corollary of the addition formula for theta functions (see e.g. \cite{BL}), we have
\begin{eqnarray*}
\theta(a_1, a_2, 1) \theta(a'_1, a_2', 1 ) & = & \theta(a_1' + a_1, a_2' + a_2, 2 ) \theta(a_1' - a_1, a_2' - a_2, 2) + \\
&& + \theta(a_1' + a_1 + 4, a_2' + a_2, 2) \theta(a_1' - a_1 + 4, a_2' - a_2, 2 ) + \\
&& + \theta( a_1'+ a_1, a_2' + a_2' + 4, 2) \theta(a_1' - a_1, a_2' - a_2 + 4, 2 ) + \\
&& + \theta(a_1' + a_1 + 4, a_2' + a_2 + 4, 2) \theta(a_1' - a_1 + 4, a_2' - a_2 + 4, 2 )\, .
\end{eqnarray*}
Using this identity and the product of Eq. (\ref{prod}), a long but straightforward computation leads to the relation $W_i = \sum_{j=0}^3 m_{ij} Z_j$, where $M:=(m_{ij})$ is the matrix
$$
M =
\left[
\begin{matrix}
8 g^3 & 4 h^2 g & 4 j^2 g & 4 k^2 g & 2 j k h \\
8 h^3 & 4 g^2 h & 4 k^2 h & 4 j^2 h & 2 j k g \\
8 j^3  & 4  k^2 j  & 4 g^2 j  & 4 h^2 j & 2 h g k \\
8 k^2 & 4 j^2 k   & 4 h^2 k & 4 g^2 k & 2 h g j 
\end{matrix}
\right] 
$$
whose entries are monomials in the following linear combination of theta functions
$$
\displaylines{
g  :=  \frac{1}{2} (\theta(0,0,2) + \theta(4,0,2) + \theta(0,4,2) + \theta(4,4,2))\, , \cr
h :=  \theta(2,0,2) + \theta (2,4,2), \quad j  :=  \theta(0,2,2) + \theta (4,2,2), \quad h :=  \theta (2,2,2) + \theta(6,2,2)\, . 
}
$$
An explicit computation shows that $\ker M$ is generated by $[1,A,B,C,D]$, where $A,B,C,D$ are as in Sec.\! \ref{kummer}. This is the expected relation in degree four. 

We conclude that the mirror of a symplectic Kummer with symplectic form (\ref{omega}) can be realized as a singular complex Kummer surface embedded in $\P^3$.

\section{Affine symplectomorphisms and algebras of Sklyanin type}\label{Sklyanin}
Consider the affine symplectomorphism
$$
\rho( x_1,  x_2 ,  y_2 ,  y_2 ) := (x_1 + b_1, x_2 + b_2, y_1 + 3x_1, y_2 + 3 x_2 )
$$
of the universal cover of the symplectic four torus $T$. \footnote{The results of this section do not extend to Kummer surfaces as translations are not compatible with the Kummer involution.} 

Under the action of $\rho$, the basic linear Lagrangian $L$ of Sec.\! \ref{Seidel} transforms into the family of affine linear Lagrangians
$$
\rho^l L =  \{ y_i = 3 l ( x - l b_i) + 3(l-1) b_i \} \qquad i=1,2 \quad l\ge1\,.
$$
If we let
$$
\yco{l}{a_1,a_2}:= (\frac{a_1}{3l} + \frac{l^2-l+1}{l} b_1, \frac{a_2}{3l} +\frac{l^2 - l + 1}{l} b_2 , 0 , 0 ) 
$$
for $a_i \in \Z$, $l \in \Z_{\ge 0}$, then the Fukaya product in $\RR(T)$ is defined by linearly extending the relations
$$
\yco{l}{a_1,a_2}\yco{l}{c_1,c_2}:=\sum_{i, j \in \Z/2\Z \times \Z/2\Z}  \theta(c_1 - a_1 + 3 i l, c_2 - a_2 + 3 j l, b_1,b_2, l) \yco{2 l}{a_1 + c_1 + 3 l i, a_2 + c_2 + 3 l j}
$$
where
$$
\theta(a_1, a_2, b_1, b_2, l) := \sum_{m,n\in Z} e^{6 l \pi i \tau_1 (m+\frac{a_1}{6 l}+\frac{b_1}{2})^2} e^{6l \pi i \tau_2(n+\frac{a_2}{6l}+\frac{b_2}{2})^2} e^{12l \pi i \tau_3(m+\frac{a_1}{6l}+\frac{b_1}{2})(n+\frac{a_2}{6l}+\frac{b_2}{2})} \, .
$$
Since $\theta(a_1,a_2,0,0,1)=\theta(6-a_1,6-a_2,0,0,1)$, the product is commutative when $b_1, b_2\in \Z$. In this case, there is a surjective ring homomorphism $\C[Z_{ij}]\to \RR(T)$ or equivalently, the abelian surface mirror to $T$ is embedded into $\P^8$. Here $\C[Z_{ij}]$ is the free commutative polynomial ring with generators $Z_{ij}:=\yco{1}{i,j}$ for $i,j\in \Z/3\Z$, which is obviously isomorphic to the quotient of the free associative ring $\C \langle Z_{ij} \rangle$ by the ideal generated by the relations $[Z_{ij}, Z_{kl}]=0$. \footnote{The trivial relations for $(i,j) = (k,l)$ are obviously not affected by noncommutative deformations.} On the other hand, when $b_1, b_2 \notin \Z$, there is a surjective homomorphism of graded associative rings $\C \langle Z_{ij}\rangle/(I)  \to \RR(T)$. To define the ideal $I$ we need some notation. 
We define vectors
\begin{eqnarray*}
v_1^\pm & : = & ( A_{20} \pm A_{40} , A_{50} \pm A_{10}, A_{23} \pm A_{43}, A_{53} \pm A_{13}) ^T \\
v_2^\pm & : = & ( A_{02} \pm A_{04} , A_{32} \pm A_{34}, A_{05} \pm A_{01}, A_{35} \pm A_{31}) ^T \\
v_3^\pm & : = & ( A_{22} \pm A_{44} , A_{52} \pm A_{14}, A_{25} \pm A_{41}, A_{55} \pm A_{11}) ^T \\
v_4^\pm & : = & ( A_{24} \pm A_{41} , A_{54} \pm A_{12}, A_{21} \pm A_{45}, A_{51} \pm A_{15}) ^T \\
v_5 & : = & ( A_{00} , A_{30} , A_{03}, A_{33}) ^T 
\end{eqnarray*} 
where $A_{ij} := \theta(i,j,b_1,b_2,1)$. Consider the matrix $V_j$ obtained from $(v_1^+|v_2^+|v_3^+|v_4^+|v_5)$ by removing the j-th column and define the $4\times 5$ matrix $\bf M$ whose $i$-th column $m_i$ satisfies $V_i m_i = - v_i^-$.  

Moreove,r let us define a vector of commutators
$$
\ZZ_{ij}^- : = 
\left(
\begin{matrix}
Z_{(i-1)j} Z_{(i+1)j} - Z_{(i+1)j} Z_{(i-1)j} \\
Z_{i(j-1)} Z_{i(j+1)} - Z_{i(j+1)} Z_{i(j-1)}\\
Z_{(i-1)(j-1)} Z_{(i+1)(j+1)} - Z_{(i+1)(j+1)} Z_{(i-1)(j-1)}\\
Z_{(i-1)(j+1)} Z_{(i+1)(j-1)} - Z_{(i+1)(j-1)} Z_{(i-1)(j+1)}
\end{matrix}
\right)
$$
and a vector of anticommutators
$$
\ZZ_{ij}^+ : = 
\left(
\begin{matrix}
Z_{(i-1)j} Z_{(i+1)j} + Z_{(i+1)j} Z_{(i-1)j} \\
Z_{i(j-1)} Z_{i(j+1)} + Z_{i(j+1)} Z_{i(j-1)}\\
Z_{(i-1)(j-1)} Z_{(i+1)(j+1)} + Z_{(i+1)(j+1)} Z_{(i-1)(j-1)}\\
Z_{(i-1)(j+1)} Z_{(i+1)(j-1)} + Z_{(i+1)(j-1)} Z_{(i-1)(j+1)}\\
Z_{ij}Z_{ij}
\end{matrix}
\right)\, .
$$
The ideal $I$ is then generated by relations $\ZZ_{ij}^- = \M \ZZ_{ij}^+$ for all $(i,j)\in \Z/3\Z \times \Z/3\Z$. For any $\tau_1$, $\tau_2$, $\tau_3$ and $b_1,b_2\notin \Z$, we interpret $\C\langle Z_{ij} \rangle / (I)$ as the twisted homogeneous coordinate ring of a noncommutative deformation of $\P^8$ because when $b_1, b_2 \in\Z$ we just recover $\C[Z_{ij}]$.

\end{document}